\documentclass[11pt]{article}

\input{epsf}

\usepackage{amsfonts,amssymb}

\newcommand{\qed}{\hbox{\rule{6pt}{6pt}}}

\newcommand{\Z}{\mathbb{Z}}
\newcommand{\R}{\mathbb{R}}

\newtheorem{theorem}{Theorem}[section]

\newtheorem{lemma}[theorem]{Lemma}
\newtheorem{proposition}[theorem]{Proposition}

\newtheorem{definition}[theorem]{Definition}

\newtheorem{example}[theorem]{Example}

\begin{document}

\title{Cocycle Knot Invariants, Quandle Extensions, \\ and
  Alexander Matrices}

\author{
J. Scott Carter\footnote{Supported in part by NSF Grant DMS \#9988107.}
\\University of South Alabama \\
Mobile, AL 36688 \\ carter@mathstat.usouthal.edu
\and 
Angela Harris
\\University of South Alabama \\
Mobile, AL 36688 \\ harris@mathstat.usouthal.edu
\and 
Marina Appiou Nikiforou
\\ University of South Florida
\\ Tampa, FL 33620  \\ mappiou@math.usf.edu
\and 
Masahico Saito\footnote{Supported in part by NSF Grant DMS \#9988101.}
\\ University of South Florida
\\ Tampa, FL 33620  \\ saito@math.usf.edu
}

\maketitle

\vspace{10mm}

\begin{abstract}
The  theory of quandle (co)homology  and cocycle knot invariants
 is rapidly being developed.  We begin with a summary of these recent 
advances. One such advance is the notion of a dynamical cocycle.
 We show how dynamical cocycles can be used to color knotted surfaces
 that are obtained from classical knots by twist-spinning. 
We also 
demonstrate relations between cocycle invariants and Alexander 
matrices. 
\end{abstract}

\vspace{10mm}

\section{Introduction}
The first half of 
this 
paper is a survey of the rapidly growing area of knot invariants 
and knotted surface invariants that are defined via quandles and their
 cocycles. 
Several key examples 
are 
closely examined from 
the  viewpoint of these 
recent developments.
In particular, 
dynamical cocycles are used to color knotted surfaces
 that are obtained from classical knots by twist-spinning,
 and   relations between cocycle invariants and Alexander
matrices are demonstrated. 
A quandle is a set with a 
self-distributive binary operation (defined below)
whose definition was partially motivated from knot theory. 
A (co)homology theory was defined in \cite{CJKLS} for quandles,
which is a modification of rack (co)homology defined in \cite{FRS}. 
The cohomology theory has found applications to the classification 
of Nichols 
algebras \cite{AG}. 
State-sum invariants,
called the quandle cocycle invariants, 
 using quandle cocycles as weights are 
defined \cite{CJKLS} and computed for important families
of classical knots and knotted surfaces \cite{CJKS1}.
Other 
survey articles on this 
subject are available \cite{CS:sur,Kam:sur}.

In this paper, first we give a short 
overview of 
the subject in Sections~\ref{qsec}, \ref{homolsec}, 
\ref{oldinvsec}, and \ref{extsec}.
New extensions called 
extensions by dynamical cocycles 
defined in \cite{AG} 
are 
studied in relation to colorings of twist-spun knots 
in Section~\ref{twistspunsec}.
Then, we define a generalized cocycle invariant in
the form of a family of vectors, 
which 
encorporates the 
refined version given  
in \cite{Lopes}, with
detailed computations  for a few examples in Section~\ref{invsec}. 
 Relations of these cocycle invariants and Alexander 
matrices are 
given in Section~\ref{LXsec}.

\section{Quandles and Quandle Colorings} \label{qsec}

In this section we 
define quandles and quandle colorings.

A {\it quandle}, $X$, is a set with a binary operation 
$(a, b) \mapsto a * b$
such that

(I) For any $a \in X$,
$a* a =a$.

(II) For any $a,b \in X$, there is a unique $c \in X$ such that 
$a= c*b$.

(III) 
For any $a,b,c \in X$, we have
$ (a*b)*c=(a*c)*(b*c). $

A {\it rack} is a set with a binary operation that satisfies 
(II) and (III).

\begin{figure}[h]
\begin{center}
\mbox{
\epsfxsize=4in
\epsfbox{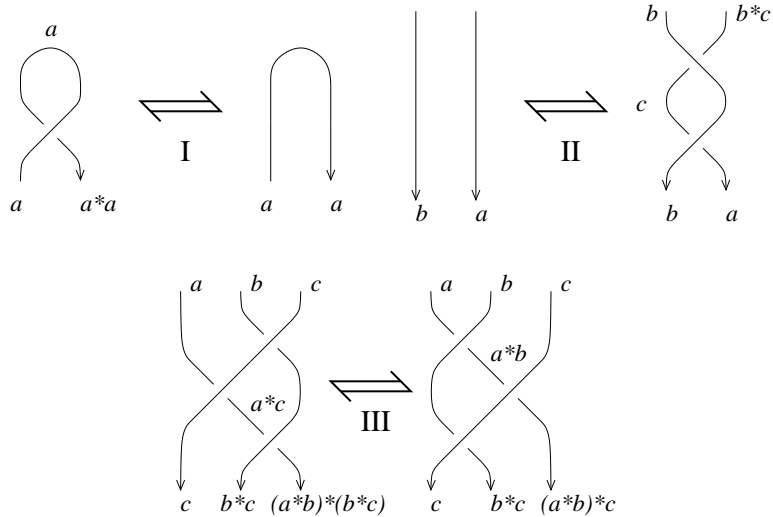}
}
\end{center}
\caption{Reidemeister moves and quandle conditions}
\label{Rmoves}
\end{figure}

Racks and quandles have been studied in, for example, 
\cite{Br88,FR,Joyce,K&P,Matveev}.
The axioms for a quandle correspond respectively to the 
Reidemeister moves of type I, II, and III 
(see Fig.~\ref{Rmoves} and 
\cite{FR,K&P}, for example). 
Quandle structures
 have been found in areas other than knot theory, 
see  \cite{AG} 
and  \cite{Br88} 
for example. Recently, it was pointed out in \cite{Yetter}
that  simple curves on a surface have a quandle structure
via an action by Dehn twists.

A function $f: X \rightarrow  Y$ between quandles
or racks  is a {\it homomorphism}
if $f(a \ast b) = f(a) * f(b)$ 
for any $a, b \in X$. 
The following are typical examples of quandles.

\begin{enumerate}   
\item
A group $X=G$ with
$n$-fold conjugation
as the quandle operation: $a*b=b^{-n} a b^n$.

\item Any subset of $G$ that is closed under conjugation.

\item
Let $n$ be a positive integer.
For elements  $i, j \in \{ 0, 1, \ldots , n-1 \}$, define
$i\ast j \equiv 2j-i \pmod{n}$.
Then $\ast$ defines a quandle
structure  called the {\it dihedral quandle},
  $R_n$.
This set can be identified with  the
set of reflections of a regular $n$-gon
  with conjugation
as the quandle operation.
\item
Any $\Lambda (={\Z }[T, T^{-1}])$-module $M$
is a quandle with
$a*b=Ta+(1-T)b$, $a,b \in M$, called an {\it  Alexander  quandle}.
Furthermore for a positive integer
$n$, a {\it mod-$n$ Alexander  quandle}
${\Z }_n[T, T^{-1}]/(h(T))$
is a quandle
for
a Laurent polynomial $h(T)$.
It 
is finite if the coefficients of the
highest and lowest degree terms
of $h$   are units in $\Z_n$.
The dihedral quandle $R_n$ can be identified with 
${\Z}_n [T,T^{-1}]/(T+1)$. 
\end{enumerate}

\begin{figure}
\begin{center}
\mbox{
\epsfxsize=2.5in
\epsfbox{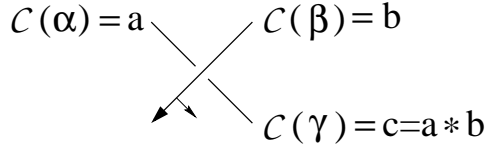} 
}
\end{center}
\caption{ Quandle relation at a crossing  }
\label{qcolor} 
\end{figure}

Let $X$ be a fixed quandle.
Let $K$ be a given oriented classical knot or link diagram,
and let ${\cal R}$ be the set of (over-)arcs. 
The normals are given in such a way that (tangent, 
normal) matches
the orientation of the plane, see Fig.~\ref{qcolor}. 
A (quandle) {\it coloring} ${\cal C}$ is a map 
${\cal C} : {\cal R} \rightarrow X$ such that at every crossing,
the relation depicted in Fig.~\ref{qcolor} holds. 
More specifically, let $\beta$ be the over-arc at a crossing,
and $\alpha$, $\gamma$ 
be under-arcs such that the normal of the over-arc
points from $\alpha$ to  $\gamma$.
Then it is required that ${\cal C}(\gamma)={\cal C}(\alpha)*{\cal C}(\beta)$.
The color ${\cal C}(\gamma)$ depends only on the choice 
of orientation of the over-arc; therefore this rule defines the coloring
 at both positive and negative crossings.

For example, 
Fox's $n$-coloring 
\cite{FoxTrip} 
is a quandle coloring by the dihedral quandle $R_n$.
The classical result that a knot is 
non-trivially 
Fox $n$-colorable 
(for $n$ prime) 
if $n| \Delta(-1)$ 
(where $\Delta(T)$  
denotes the Alexander polynomial)
has been generalized by Inoue~\cite{Inoue} to the following:

Let $\Delta_K^{(i)}(T)$ denote the greatest common divisor of 
all $(n-i-1)$ minor determinants of the presentation matrix 
for the knot module obtained via the Fox calculus.

\begin{theorem}{\bf \cite{Inoue}} \label{Inouethm}
Let $p$ be a prime number, 
$J$ an ideal of the ring $\Lambda_p=\Z_p[T,T^{-1}]$.  
For each 
$i\ge 0$, put  
$e_i(T)=\Delta_{K}^{(i)}(T)/\Delta_{K}^{(i+1)}(T).$
Then the number of 
colorings
by 
the Alexander quandle 
$\Lambda_p/J$ is equal to the cardinality of the module 
$\Lambda_p/J \oplus \oplus_{i=0}^{n-2} \{ \Lambda_p/(e_i(T),J)\}$.
\end{theorem}

Alternatively, 
a coloring can be described as a 
quandle homomorphism as follows.
Classical knots have fundamental quandles that are defined via generators
and relations. The theory of quandle presentations is given a complete
treatment in \cite{FR}. 
The quandle relation
$a*b=c$  holds where $a$ is the generator that corresponds to the 
under-arc 
away from which the normal to the over-arc  
points, $b$ is the generator
that corresponds to the over-arc
and $c$ corresponds to the under-arc 
towards which the transversal's normal points,
see Fig.~\ref{qcolor}.
A  coloring  of a classical knot diagram by 
a quandle $X$ gives rise to 
a quandle homomorphism from the fundamental quandle to the quandle $X$.

Using Waldhausen's theorem Joyce shows:
\begin{theorem}{\bf \cite{Joyce,Matveev}} 
If two knots in $\R^3$ 
have isomorphic 
fundamental 
quandles, then the 
knots are equivalent up to orientations of $\R^3$ and the knots. 
\end{theorem}
This fundamental fact was generalized:
\begin{theorem}
{\bf \cite{FR}} The fundamental augmented rack is a complete
 invariant for irreducible semi-framed links in closed connected $3$-manifolds.
\end{theorem}
See \cite{FR} for a full account of the notation and terminology.
Using an interpretation of cocycle
knot invariants in terms  of 
the {\it canonical} class $c(L)$ of a link $L$, 
the above  theorem was further generalized to:
\begin{theorem}{\bf \cite{RS}}
If $L,M$ are two links in $S^3$ such that there is an isomorphism  $\phi$ of 
fundamental racks with $\phi_* (c(L))=c(M)$, then 
$L$ and $M$ are isotopic.
\end{theorem}

\section{Quandle Homology and Cohomology Theories} \label{homolsec}

In this section, we present  
the  ordinary 
 quandle homology theory.
Originally, rack homology and homotopy theory were defined and 
studied in \cite{FRS}, and a modification to quandle homology theory 
was given in \cite{CJKLS} to define a knot invariant in a state-sum form.
Then they were  generalized to 
a  twisted theory in \cite{CENS}. 
The most general form of the quandle homology known to date is given
 in \cite{AG}.
Computations are found in \cite{CJKS1,betti}  and also in
\cite{EG,LN,Mochi} by other authors.

 Let $C_n^{\rm R}(X)$ be the free 
abelian group generated by
$n$-tuples $(x_1, \dots, x_n)$ of elements of a quandle $X$. Define a
homomorphism
$\partial_{n}: C_{n}^{\rm R}(X) \to C_{n-1}^{\rm R}(X)$ by \begin{eqnarray}
\lefteqn{
\partial_{n}(x_1, x_2, \dots, x_n) } \nonumber \\ && =
\sum_{i=2}^{n} (-1)^{i}\left[ (x_1, x_2, \dots, x_{i-1}, x_{i+1},\dots, x_n) \right.
\nonumber \\
&&
- \left. (x_1 \ast x_i, x_2 \ast x_i, \dots, x_{i-1}\ast x_i, x_{i+1}, \dots, x_n) \right]
\end{eqnarray}
for $n \geq 2$ 
and $\partial_n=0$ for 
$n \leq 1$. 
 Then
$C_\ast^{\rm R}(X)
= \{C_n^{\rm R}(X), \partial_n \}$ is a chain complex.

Let $C_n^{\rm D}(X)$ be the subset of $C_n^{\rm R}(X)$ generated
by $n$-tuples $(x_1, \dots, x_n)$
with $x_{i}=x_{i+1}$ for some $i \in \{1, \dots,n-1\}$ if $n \geq 2$;
otherwise let $C_n^{\rm D}(X)=0$. If $X$ is a quandle, then
$\partial_n(C_n^{\rm D}(X)) \subset C_{n-1}^{\rm D}(X)$ and
$C_\ast^{\rm D}(X) = \{ C_n^{\rm D}(X), \partial_n \}$ is a sub-complex of
$C_\ast^{\rm
R}(X)$. Put $C_n^{\rm Q}(X) = C_n^{\rm R}(X)/ C_n^{\rm D}(X)$ and 
$C_\ast^{\rm Q}(X) = \{ C_n^{\rm Q}(X), \partial'_n \}$,
where $\partial'_n$ is the induced homomorphism.
Henceforth, all boundary maps will be denoted by $\partial_n$.

For an abelian group $G$, define the chain and cochain complexes
\begin{eqnarray}
C_\ast^{\rm W}(X;G) = C_\ast^{\rm W}(X) \otimes G, \quad && \partial =
\partial \otimes {\rm id}; \\ C^\ast_{\rm W}(X;G) = {\rm Hom}(C_\ast^{\rm
W}(X), G), \quad
&& \delta= {\rm Hom}(\partial, {\rm id})
\end{eqnarray}
in the usual way, where ${\rm W}$ 
 $={\rm D}$, ${\rm R}$, ${\rm Q}$.

The groups of cycles and boundaries are denoted respectively by 
${\mbox{\rm ker}}(\partial) =Z_n^{\rm W}(X;G) \subset C_n^{\rm W} (X;G)$ and 
${\mbox{\rm Im}}(\partial)=B_n^{\rm W}(X;G) \subset C_n^{\rm W}(X;G)$ 
while the cocycles and coboundaries are denoted respectively by
${\mbox{\rm ker}}(\delta) =Z^n_{\rm W}(X;G) \subset C^n_{\rm W}(X;G)$ and 
${\mbox{\rm Im}}(\partial)=B^n_{\rm W}(X;G) \subset C^n_{\rm W}(X;G).$
 In particular,  a quandle $2$-cocycle is
an element  $\phi \in Z^2_{\rm Q}(X;G)$,  and 
the equalities
\begin{eqnarray*}
\phi(x,z)+\phi(x*z,y*z)&=&\phi(x*y,z)+\phi(x,y)  \\
 \mbox{and} \quad  \phi(x,x) & = & 0
\end{eqnarray*}
are satisfied for all $x,y,z\in X$.

The $n$\/th {\it quandle homology group\/}  and the $n$\/th
{\it quandle cohomology group} \cite{CJKLS} of a quandle $X$ with coefficient group $G$ are
\begin{eqnarray}
H_n^{\rm Q}(X; G) 
 &= &H_{n}(C_\ast^{\rm Q}(X;G)) \; = \; Z_n^{\rm Q}(X;G)/B_n^{\rm Q}(X;G), 
 \nonumber \\
H^n_{\rm Q}(X; G) 
 &= & H^{n}(C^\ast_{\rm Q}(X;G))\; = \; Z^n_{\rm Q}(X;G)/B^n_{\rm Q}(X;G). \end{eqnarray}

The following developments have been made recently.

\begin{itemize}
\item  
The conjecture made in \cite{betti} on  
the long 
exact homology 
sequence
was proved 
by Litherland and Nelson 
in \cite{LN}: 
 the long exact sequence  of quandle homology 
$$\cdots 
\rightarrow H^{\rm D}_n(X;A) \rightarrow H^{\rm R}_n(X;A)
\rightarrow H^{\rm Q}_n(X;A) \rightarrow  H^{\rm D}_{n-1}(X;A)\rightarrow \cdots $$
splits
into short exact sequences 
$$ 0 {\to}  H_n^{\rm D}(X;A) {\to}
H_n^{\rm R}(X;A)
{\to} H_n^{\rm Q}(X;A)
{\to} 0 . $$

\item
Mochizuki has computed several key cohomology groups. 
We highlight some of his results.
First,
$$H^3_{\rm Q}(R_p;\Z_p) \cong 
\Z_p,$$  
and he gives an explicit expression for a generating cocycle. 
He gives explicit $2$-cocycles for Alexander quandles over a 
field $K$ thereby computing
$H^2_{\rm Q}( K[T,T^{-1}]/(T-\omega);K)$ for $\omega \ne 0,1.$
He shows that dihedral quandles of odd order have vanishing rational
 cohomology in all dimensions. This was shown independently in \cite{LN}.

\item Etingof and Gra\~{n}a \cite{EG} have computed the rational
 cohomology of any finite quandle in all dimensions. In particular, 
they show that the bounds on the rank of the betti numbers given in
 \cite{betti} are equalities.
Furthermore, they relate the $2$-dimensional quandle homology to group 
cohomology. Let $X$ be a quandle and $G_X$ be its enveloping group:
 $G_X = \langle x\in X: y^{-1}x y = x*y \rangle.$ 
If $A$ is a trivial $G_X$-module (as is the case with ordinary quandle
 homology), then
$$H^2_{\rm Q}(X;A) \cong 
 H^1 (G_X; \mbox{Fun}(X,A) )$$ 
 where  
$\mbox{Fun}(X,A)$
denotes 
the set of functions.

\item Andruskiewitsch and Gra\~{n}a  \cite{AG} have developed the
 theory of quandle and rack cohomology further.
 They have developed a cohomology theory that 
encompasses those 
in \cite{CENS} and Ohtsuki's theory \cite{Ohtsuki}. 
Furthermore, their primary interest is in the classification of
 certain pointed Hopf
 algebras called Nichols algebras.
The quandle cohomology plays a central role here.

\item 
Using dynamical  cocycles, Gra\~{n}a \cite{G2} 
classified indecomposable racks of order $p^2$ for any prime $p$.
Another classification theorem was proved by Nelson \cite{Nel}
who described isomorphism classes of Alexander quandles 
by the submodules $\mbox{Im}(1-T)$.

\end{itemize}

\section{Cocycle Knot Invariants} \label{oldinvsec}
\label{quanss}

\subsection*{Cocycle knot invariants of classical and virtual knots}

Let $K$ be a classical knot or link diagram. Let a finite 
quandle $X$, and an (untwisted) 
quandle $2$-cocycle $\phi \in Z^2_{\rm Q}(X;A)$ be given.
A {\it (Boltzmann) weight}, $B(\tau, {\cal C})$ 
(that depends on $\phi$),
at a  crossing $\tau$ is defined as follows.
Let  ${\cal  C}$ 
denote a coloring ${\cal  C}: {\cal R} \rightarrow X$. 
Let $\beta$ be the over-arc at $\tau$, and $\alpha$, $\gamma$ be 
under-arcs such that the normal to $\beta$ points from $\alpha$ to $\gamma$,
see Fig.~\ref{qcolor}. 
Let $x={\cal C}(\alpha)$ and $y={\cal C}(\beta)$. 
Then define $B(\tau, {\cal C})= \phi(x,y)^{\epsilon (\tau)}$,
where $\epsilon (\tau)= 1$ or $-1$, if  
(the sign of) the crossing  $\tau$ 
is positive or negative, respectively.
By convention, the crossing in Fig.~\ref{qcolor} is positive if 
the orientation of the under-arc points downward.

\begin{figure}
\begin{center}
\mbox{
\epsfxsize=4in
\epsfbox{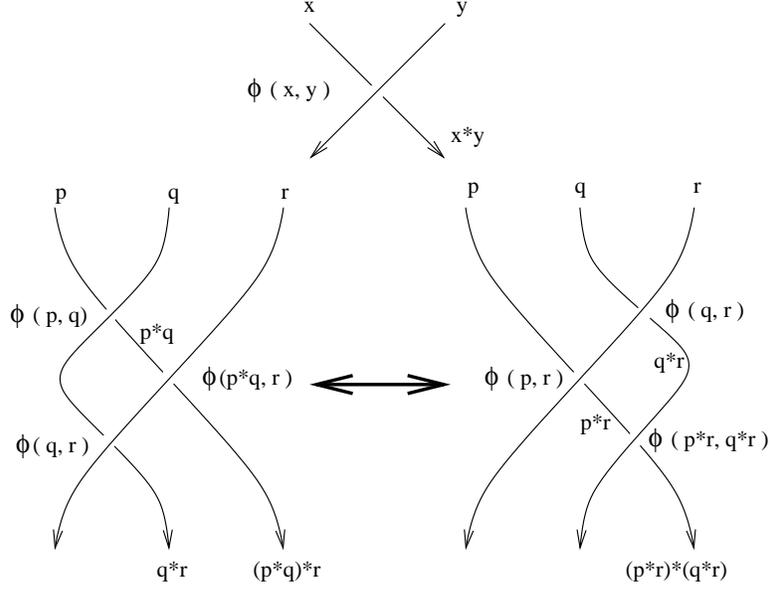} 
}
\end{center}
\caption{ The untwisted $2$-cocycle condition and type III move  }
\label{2cocy} 
\end{figure}

The {\it (quandle) cocycle knot invariant} is defined by 
the state-sum expression  
$$
\Phi (K) = \sum_{{\cal C}}  \prod_{\tau}  B( \tau, {\cal C}). 
$$
The product is taken over all crossings of the given diagram $K$,
and the sum is taken over all possible colorings.
The values of the partition function 
are  taken to be in  the group ring ${\Z}[A]$ where $A$ is the coefficient 
group  written multiplicatively. 
The state-sum depends on the choice of $2$-cocycle $\phi$. 
This is proved in
\cite{CJKLS} to be a knot invariant.  
Figure~\ref{2cocy} shows the invariance of the state-sum under
the Reidemeister type III move.
The sums of cocycles, equated before and  after the move, 
is the $2$-cocycle condition.

Relations to braid group representations and 
quantum invariants are studied in \cite{Grana}, 
see also \cite{CS:sur} for a viewpoint 
from  the bracket state-sum form and Dijkgraaf-Witten invariants.

\subsection*{Cocycle 
Invariants for Knotted Surfaces}

The state-sum invariant is defined in an analogous way for 
oriented knotted surfaces
in $4$-space using their projections and diagrams in $3$-space. 
Specifically, the above steps can be repeated as follows,
for a fixed finite quandle $X$ and a knotted surface diagram $K$.

\begin{figure}[h]
\begin{center}
\mbox{
\epsfxsize=3in
\epsfbox{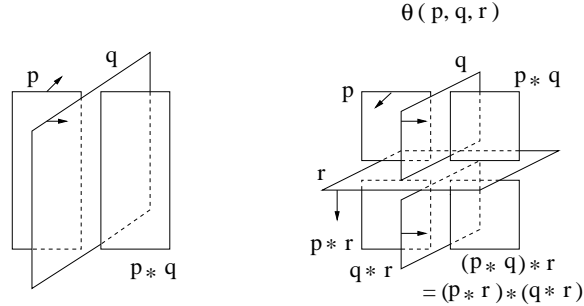} 
}
\end{center}
\caption{Colors at double curves and $3$-cocycle at a triple point }
\label{triplept} 
\end{figure}

\begin{itemize}
\item
The diagrams consist of double curves and isolated branch and triple points
\cite{CS:book}. Along the double curves, 
the 
coloring rule is defined using normals
in the same way as classical case, as depicted in
the left of Fig.~\ref{triplept}.

\item
The sign $\epsilon(\tau)$ of a triple point $\tau$
is  defined \cite{CS:book} 
in such a way that it is positive if and only if 
the normals to top, middle, bottom sheets, in this order, 
match the orientation of $3$-space. 

\item  
For a coloring ${\cal C}$, the Boltzman weight 
at a triple point $\tau$ is defined by 
$B(\tau, {\cal C})=$ \linebreak 
$\theta(x,y,z)^{\epsilon (\tau)} $, 
where $\theta$ is a $3$-cocycle, $\theta \in Z^3_{\rm Q}(X;A)$. 
In the right of Fig.~\ref{triplept},
the triple point $\tau$ is positive, 
so that $B(\tau, {\cal C})=\theta(p,q,r)$. 

\item 
The state-sum is defined by $\Phi(K)=  
\sum_{{\cal C}}  \prod_{\tau}  B( \tau, {\cal C}). 
$

\end{itemize}

Recall that a function $\theta:X\times X\times X \rightarrow A$
 is a quandle $3$-cocycle 
if  
\begin{eqnarray*}
\theta (p,r,s)  +\theta(p*r,q*r,s) +\theta(p,q,r) 
  &=& \theta(p*q,r, s) + \theta(p,q,s) + \theta(p*s,q*s,r*s) \\
\theta(p,p,q) & = & 0 \\
\theta(p,q,q) &=& 0 \end{eqnarray*}
By checking the analogues of Reidemeister moves for knotted
surface diagrams, called Roseman moves, 
it was shown in \cite{CJKLS} that 
$\Phi(K)$ is an invariant, 
called  the {\it (quandle) cocycle invariant}
of knotted surfaces. 

The value of the state-sum invariant depends
only on the cohomology class represented by the defining cocycle.
In particular, a coboundary will simply count the number of colorings
of a knot 
or knotted surface 
by the quandle $X$.

\subsection*{Applications}

Important topological applications have been obtained
using the cocycle invariants for knotted surfaces.

\begin{itemize}
\item
The $2$-twist spun trefoil $K$ and its orientation-reversed counterpart $-K$
have shown to have distinct cocycle invariants
using a cocycle in $Z^3_Q(R_3;\Z_3)$, 
providing a proof that $K$ is non-invertible \cite{CJKLS}.
The higher genus surfaces obtained from $K$ by adding arbitrary
number of trivial $1$-handles are also non-invertible, since 
such handle additions do not alter the cocycle invariant.
This result in higher genus cases
is not immediately obtained from \cite{Hillman,Ruber}, although
higher genus generalizations of the Farber-Levine pairing
\cite{Kawa90a} can be used. 

\item 
Cocycle invariants for twist spun $(2,n)$-torus knots
were computed using Maple \cite{CJKS1} for some quandles.
Computer-free calculations and general formulas were 
obtained later in \cite{Satoh} using 
explicit formulas of $3$-cocycles provided in \cite{Mochi}
for dihedral quandles.  
Mochizuki's formulas were also used for  the following geometric application.  

\item
The projection of the $2$-twist spun trefoil was shown to 
have at least four triple points \cite{SatShima}.

\item 
The projection of the $3$-twist spun trefoil was shown to 
have at least six triple points \cite{SatShima2}. 
The cocycle employed here appeared in \cite{CJKS1}.
 Satoh and Shima 
gave a set of linear equations among numbers of colored
triple points to give 
algebraic lower bounds on the number of triple points.
 They have developed a computer program to compute these bounds.

Colorings of knotted surfaces in relation to dynamical cocycles
are discussed in Section~\ref{twistspunsec}, 
that would provide the first step towards 
extending their results. 

\end{itemize}

\section{ Extension theory of quandles} \label{extsec}

Let $X$ be a quandle, and for a given abelian coefficient group $A$, take 
a $2$-cocycle $\phi \in Z^2_{\rm Q}(X;A)$. 
Let $E=A \times X$ and define a binary operation by 
$(a_1, x_1)*(a_2, x_2)=(a_1 + \phi(x_1, x_2), x_1 * x_2)$.
It was shown in \cite{CENS} that $(E, *)$ defines a quandle, called
 an {\it abelian (\mbox{\rm or} central) extension}, and is 
denoted by $E=E(X, A, \phi) $. 
(This  is in parallel to central extension of groups,
see Chapter IV of \cite{Brown}.)
The following examples were given in \cite{CENS}.

\begin{itemize}

\item 
For any positive integer $q$ and $m$, the quandle 
$E=W_{m+1}=\Z_q [T, T^{-1} ] / (1-T)^{m+1} $ is an abelian extension
of $X=W_{m}=\Z_q [T, T^{-1} ] / (1-T)^{m} $ over $\Z_q$: 
$E=E(X, \Z_q, \phi)$, for some $\phi \in Z^2_{\rm Q}(X; \Z_q)$.

\item
For any  positive integers $q$ and $m$, 
$E=U_{m+1}=\Z _{q^{m+1}} [T, T^{-1}] /  (T -1 +q) $ is an abelian extension 
$E=E( \Z _{q^{m}} [T, T^{-1}] /  (T -1 +q) ,  \Z_q,  \phi)$
of $X=U_{m}= \Z _{q^{m}} [T, T^{-1}] /  (T -1 +q)$ for some cocycle
$\phi \in  Z^2_{\rm Q}( X;  \Z_q)$. 

\end{itemize}

\noindent
For these quandles, explicit 
formulas were obtained in \cite{CENS} 
using extensions as follows. 

\begin{itemize}
\item
Represent elements of $E=W_{m+1}$ by 
$A=A_m (1-T)^m + \cdots + A_1(1- T) + A_0$, where
$A_j \in \Z_q $, $j=0, \ldots, m$.
Define  $f: E=W_{m+1} \rightarrow \Z_q \times 
X (=W_m)$ by $$f(A)=(A_m \  (\mbox{mod}\ q), 
 \  \overline{A} \  (\mbox{mod}\ (1-T)^{m-1}) ), $$
where 
$ \overline{A}=  \sum_{j=0}^{m-1} A_{j} (1-T)^{j} . $
Then for $A, B  \in E$, the quandle operation is computed 
 by
\begin{eqnarray*} 
 A*B & = & TA + (1-T)B \\
 & = &  ( A_{m}- A_{m-1}+  B_{m-1}  ) (1-T)^m +
 \sum_{j=0}^{m-1} (A_{j} - A_{j-1}+B_{j-1}) (1-T)^j , 
\end{eqnarray*} 
where $A_{-1}, B_{-1} $ are understood to be zeros in the last summation,
and the coefficients are in $\Z_q$.
In $\Z_q[T, T^{-1}]/{(1-T)^{m}}$, and we have 
$f ( A*B ) =( \phi(\overline{A}, \overline{B}), \overline{A}*\overline{B})$
 where $\phi(\overline{A}, \overline{B})= B_{m-1}- A_{m-1} $.
Hence $f$ yields an isomorphism.

\item
The cocycle $\phi$ has a description using a section. Let 
$$s: \Z_{q} [ T, T^{-1}] / (1-T)^{m} \rightarrow \Z_{q} [ T, T^{-1}] / (1-T)^{m+1}$$ be a set-theoretic section defined by 
$$s\left( \sum_{j=0}^{m-1} A_j (1-T)^j 
\right)
=  \sum_{j=0}^{m-1} A_j (1-T)^j \quad  \mbox{mod}\ (1-T)^{m+1}. $$
Then we have 
$\overline{ s(X)} * \overline{ s(Y)} =\overline{ s(X *Y)}$ for any
$X, Y \in   \Z_{q} [ T, T^{-1}] / (1-T)^{m}$,  
so that 
$[  s(X) *  s(Y)- s(X *Y)]$ is divisible by $(1-T)^m$, 
and we have 
$$ \phi(\overline{A}, \overline{B})
=  [  s(A) *  s(B)- s(A *B)]/ (1-T)^m \in \Z_q. $$

\item
For $E=U_{m+1}$, 
represent elements of $\Z _{q^{m+1}}$ by $\{ 0, 1, \ldots, q^{m+1} -1 \}$
and express them 
in their
 $q^{m+1}$-ary expansion:
$$A = A_{m} q^m + \cdots + A_1 q + A_0 \in \Z _{q^{m+1}}, $$
where $0 \leq A_j < q$, $j=0, \ldots, m$. 
Then $E=U_{m+1}$ and $X=U_m$ have a similar description 
as above.

\end{itemize}

An extension theory of quandles for ``twisted''
 cohomology cocycles was  
developed in \cite{CES}, and it provided 
  more general extension theories.
In the twisted case, the coefficient group is taken to be a
$\Lambda$-module, thus has an Alexander quandle structure, and 
the extension  $AE(X,A, \phi) = (A \times X, *)$
is defined by 
$(a_1, x_1) * (a_2, x_2) = (a_1 * a_2 + \phi(x_1, x_2), x_1 * x_2)$
 for  $\phi \in Z^2_{\rm TQ} (X; A)$,
and is called 
 an {\em Alexander extension} 
of $X$ by $(A, \phi)$,
where the subscript TQ represents the twisted theory. 
For example, $R_{p^m} $ is an Alexander extension
of $R_{p^{m-1}}$ by $R_p$: $R_{p^m}= AE(R_{p^{m-1}}, \  R_p, \ \phi)$, 
for some $\phi \in Z^2_{\rm TQ}(R_{p^{m-1}}; R_p)$.

Ohtsuki \cite{Ohtsuki}
defined a new  cohomology theory for quandles and 
an 
extension theory, 
together with a list of problems in the subject. 
Further generalizations of extensions by dynamical cocycles as 
defined in \cite{AG} will be discussed and used for coloring
twist spun knots in the next section.

\section{Extensions of quandles and colorings of twist-spun knots}
\label{twistspunsec}

\subsection*{Extensions by dynamical cocycles}

This subsection is a brief summary of a  quandle extension
theory by Andruskiewitch and Gra\~{n}a \cite{AG}. 
The notation
has 
been changed below
from that given in 
\cite{AG} 
to match our conventions in this paper.
Let $X$ be a quandle and $S$ be a non-empty set. 
Let $\alpha: X \times X \rightarrow 
\mbox{\rm Fun}(S \times S, S)=S^{S \times S}$ be a function,
so that for $\sigma, \tau \in X$ and $a, b \in S$ we have 
$\alpha_{\sigma, \tau} (s, t ) \in S$. 

Then it is checked by computations 
that $S \times X$ is a quandle by the operation
$(a, \sigma)*(b, \tau)=(\alpha_{\sigma, \tau} (a,b ) ,\sigma * \tau )$,
where $ \sigma * \tau$ denotes the quandle operation in $X$, 
if and only if $\alpha$ satisfies the following conditions:

\begin{enumerate}
\item $\quad \alpha_{\sigma, \sigma}(a,a)= a$ for all $\sigma \in X$ and 
$a \in S$; 
\item $\quad \alpha_{\sigma,\tau}(-,b): S \rightarrow S$ is a bijection for 
all $\sigma, \tau \in X$ and for all $b \in S$;
\item $\quad 
\alpha_{\sigma * \tau, \eta}(\alpha_{\sigma , \tau}(a, b), c)
=\alpha_{\sigma * \eta, \tau*\eta}( \alpha_{\sigma,\eta}(a,c),
\alpha_{\tau,\eta}(b,c) )
$
for all $\sigma, \tau, \eta \in X$ and $a,b,c \in S$. 

\end{enumerate}

Such a function $\alpha$ is called a {\it dynamical quandle cocycle}.
The quandle constructed above is denoted by $S \times_{\alpha} X$, 
and is called the {\it extension} of $X$ by a dynamical cocycle $\alpha$.
The construction is general, as they show:

\begin{lemma}{\bf \cite{AG}} \label{AGlemma} 
Let $p: Y \rightarrow X$ be a surjective quandle homomorphism 
such that the cardinality of $p^{-1}(x)$ is a constant for all $x \in X$. 
Then $Y$ is isomorphic to an extension  $S \times_{\alpha} X$ of $X$ 
by some dynamical cocycle on a set $S$.
\end{lemma}

\noindent 
\subsection*{Quandle extensions in  wreath products}

Let 
$$0\rightarrow N \stackrel{i}{\rightarrow} G \stackrel{\pi}{\rightarrow} H 
\rightarrow 1$$
be a split 
short exact sequence of groups that expresses the finite group $G$ as a 
semi-direct product $G= N \rtimes H$, so that 
we have a homomorphism $s:H\rightarrow G $ with $\pi \circ s = 1_H$.
The elements of $G$ can be written as 
pairs $(x, \sigma)$ 
 where $x\in N$ and $\sigma \in H$.
The multiplication 
rule in $G$ is given by 
$( x,  \sigma)\cdot(y, \tau)=(x\sigma(y),\sigma\tau),$ 
where $\sigma(y)$ denotes the action of $H$ on $N$ that gives $G$ the 
structure of a semi-direct product.

Let $Q$ denote a subquandle of ${\mbox{\rm Conj}}(H)$
(the group $H$ with the quandle structure given by 
conjugation). 
Thus $Q$ is a subset 
of $H$ that is closed under conjugation. Let
$\widetilde{Q}= \{ (x, \sigma) : 
\sigma \in Q \}$, then $\widetilde{Q}$ is a 
quandle by conjugation in $G$. 
The group homomorphism $\pi: G \rightarrow H$ induces the quandle 
homomorphism (denoted by the same letter) $\pi:\widetilde{Q} \rightarrow Q$. 
Lemma~\ref{AGlemma} implies

\begin{lemma} \label{wreathlemma}
Suppose the cardinality $|\pi^{-1}(\sigma)|$ is 
independent of $\sigma \in Q$. Then $\widetilde{Q}$ is an extension 
$Q$ by a dynamical cocycle 
$\alpha: Q \times Q \rightarrow S^{S \times S}$
where $S$ is a set with cardinality $|\pi^{-1}(\sigma)|$. 
\end{lemma}

Now we specialize to the case that 
 $Q$ is a subquandle of 
${\mbox{\rm Conj}}(\Sigma_n)$, 
where $\Sigma_n$ denotes the symmetric group on $n$ letters, 
and 
$N=(\Z _v)^n$ for some $v\in \{ 0, 1, \ldots \}$.
(In case $v=0$, then $N$ is the direct product of the integers, and when 
$v=1$ then $N$ is trivial.)  
The  action of $\Sigma_n$ is given by 
permutation of the factors
$\sigma(x_1,\ldots, x_n)= \sigma(\vec{x}) = 
(x_{\sigma (1)}, \ldots, x_{\sigma (n)})$, 
for $\sigma \in \Sigma_n$ and
 $\vec{x}=(x_j)_{j=1}^n \in (\Z_v )^n $. 
In this situation,
 $G=(\Z _v)^n \rtimes \Sigma_n$ is also called
a wreath product and denoted by 
$G= (\Z _v) \wr \Sigma_n$. 
Hence the group operation in $G$ is written by 
 $$( \vec{x}, \sigma) \cdot (\vec{y}, \tau) = (\vec{x}\sigma(\vec{y}), \sigma \tau)$$ where 
$\sigma(\vec{y}) = (y_{\sigma(1)},\ldots, y_{\sigma(n)}).$

It is well known that 
 elements of $G$ can be represented 
by matrices with entries in $\{x^j\}$ as follows.
First,  
we represent the cyclic group $ (\Z_v)$ multiplicatively 
as $\langle x | x^v=1 \rangle$, 
and represent 
$\vec{x} \in (\Z _v)^n$ 
by $(x^{i_1}, \ldots , x^{i_n} ).$ 
Represent $\sigma \in \Sigma_n$ by 
an $(n \times n)$-matrix $M(\sigma)$
acting on vectors of $n$ letters from the 
left. 
It is a matrix with exactly one nonzero entry in each row and  column,
and each non-zero entry is $1$.  
The pair $( \vec{x}, \sigma)$ 
is  represented as a matrix  $M(\vec{x}, \sigma)$ 
obtained from  $M(\sigma)$ by replacing the non-zero entry in the $j\/$th row 
by  $x^{i_j}$. 
The group composition in $G$ is matrix multiplication of $M(\vec{x}, \sigma)$s
where $0\cdot x^i =0$ and $x^ix^j=x^{i+j}$. 
For example, we write
$$((x^i,x^j,x^k),(12))
=\left( \begin{array}
{ccc}0&x^i&0 \\x^j&0&0\\0&0&x^k 
\end{array}\right),$$
and
$$((x^\ell,x^m,x^p),(123))
=\left( \begin{array}{ccc}0&0&x^\ell \\x^m&0&0\\0&x^p&0 
\end{array}\right).$$
The matrix product evaluates as 
$$\left( 
\begin{array}{ccc}0&x^i&0 \\x^j&0&0\\0&0&x^k 
\end{array}\right)\cdot\left(
 \begin{array}{ccc}0&0&x^\ell \\x^m&0&0\\0&x^p&0 
\end{array}\right)=\left( 
\begin{array}{ccc}x^{i+m}&0&0\\0&0&x^{\ell+j}\\0&x^{k+p}&0 
\end{array}\right)$$
which 
corresponds to $((x^{i+m},x^{\ell+j},x^{k+p}),(23)).$
Meanwhile, 
the product 
$$((x^\ell,x^m,x^p),(123))\cdot((x^i,x^j,x^k),(12))= 
((x^{\ell+k},x^{m+i},x^{p+j}),(13))$$
 is represented by the matrix
$$\left(\begin{array}{ccc}0&0&x^{\ell+k} \\0&x^{m+i}&0\\x^{p+j}&0&0 \end{array}\right).$$
Note that the cycles of $\Sigma_n$ act 
on $n$-letters from the left in this convention.

We take our subquandle $Q$ to be, 
for example, the dihedral quandle, $R_n$, 
or a subset of a given conjugacy class that is itself closed under 
conjugation, 
such as  
$QS_4 = \{ (123),(142),(134),(243)  \}\subset  \Sigma_4$.
We are interested in applications to knots herein, so we 
assume that $n$ is odd.
Then for both of 
$Q=R_n$ or $QS_4$, an element $\sigma \in Q$ has a 
fixed point in $\{1,\ldots, n\}$, and the matrix representation
 $M(\sigma, \vec{x})$ 
of  $(\sigma, \vec{x})$ has exactly one element along the diagonal. 
It is easy to 
see that the exponent of the diagonal element is fixed under the conjugation 
action, so we 
restrict our attention to the subquandle $Q(v)$
of $\widetilde{Q}$ 
in which the non-zero diagonal element is  $1$.  
By Lemma~\ref{AGlemma}, we  see that  $Q(v)$ 
is also an extension of $Q$ by a dynamical cocycle.

In the case of the dihedral quandle $Q=R_n$ for $n$ odd, we 
simplify the notation further. 
Consider a regular $n$-gon whose vertices labeled with $\{1, \ldots, n\}$ 
in this order, on which $R_n$ acts as reflections. 
For $i=1, \ldots, n$, let 
$\sigma_i$ denote the reflection 
of a regular $n$-gon 
 which fixes 
the vertex 
labeled $i$. 
Then 
$ \sigma_i * \sigma_j=\sigma_{2j -i}$. 
Denote the element $(\vec{x}, \sigma_j) \in G$ 
where 
$\vec{x}=(x^{i_1}, \ldots, x^{i_{j-1}}, 1,  x^{i_{j+1}}, \ldots,x^{i_n})$,
by $\sigma_j(i_1,\ldots, \widehat{\imath}_j, \ldots, i_n)$
where $\widehat{\imath_j}$ indicates that the  
$j\/$th element in this list is missing.
For example, set $a=\sigma_1$, $b=\sigma_2$, and $c=\sigma_3$
in $R_3$ and the elements of  $R_3(v)$
are denoted by $a_{j,k}$, $b_{i,k}$, and $c_{i,j}$,
where $i,j,k \in \Z/v$. 
For convenience we summarize 
the multiplication table for  $R_3(v)$ as follows, where $r*c$ indicates that 
the table represents 
(row)$*$(column).

\vspace{2mm}
\begin{center}
{\begin{tabular}{|l||c|c|c||}\hline \hline
$r*c$ & $a_{n,\  p}$ & $b_{m,\ p}$ & $c_{m,\ n}$  \\ \hline \hline 
$a_{j,\ k}$ & $a_{k+n-p,\ j-n+p}$ & $c_{k-p,\ j+p}$ & $b_{j-n, \ k+n}$ \\ \hline
$b_{i,\ k}$ &$c_{i+p, \ k-p}$ & $b_{k+m-p,\ i-m+p}$ & $a_{i-m, \ k+m}$ \\ \hline
$c_{i,\ j}$ & $b_{i+n,\ j-n}$  & $a_{j+m,\ i-m } $ 
&  $c_{j+m-n,\ i-m+n}$ \\
\hline \hline
\end{tabular}}
\end{center}
\vspace{2mm}

\subsection*{Coloring twist-spun knots by extended quandles}

Now we use the above extensions of quandles to color twist-spun knots.

\begin{figure}[h]
\begin{center}
\mbox{
\epsfxsize=4in
\epsfbox{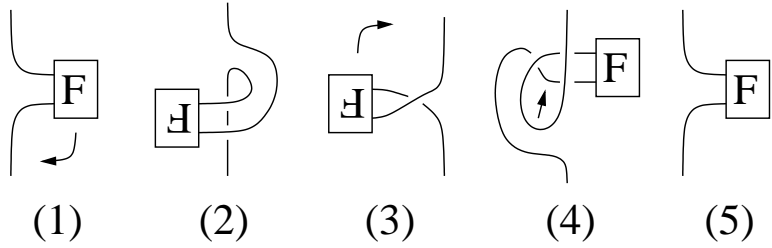} 
}
\end{center}
\caption{The general method of twist spinning}
\label{twist} 
\end{figure}

\begin{example} 
\begin{enumerate}
\item %%%1
 The $2v$-twist spun trefoil is non-trivially 
 colorable by the quandle
$R_3(v)$.

\vspace{-3mm}

\item %%%%2
The $3u$-twist spun trefoil is non-trivially colorable by
 the quandle $QS_4(u).$

\vspace{-3mm}

\item  %%%%3
The $2v$-twist spun figure 8 knot is non-trivially
 colorable by $R_5(v)$.

\end{enumerate}
 \end{example}
{\it Proof.} The schematic diagram indicated in Fig.~\ref{twist} illustrates 
Satoh's method \cite{Satoh:Korea} for obtaining the
twist spun knot 
from a 
$(1-1)$-tangle $F$ 
whose closure is 
a given classical knot 
$K$. 

{}From left to right in the diagram, a movie of one full twist 
of a tangle is depicted. After repeating $k$ full twists,
 the tangle is identified with the original one to form the $k$-twist
 spun knot of $K$. Reidemeister moves performed in the course
of the isotopy 
correspond
to 
critical or singular points on 
the projection, see \cite{CS:book}, for example,  for details.
In particular, triple points on the projection correspond to
type III moves, and they appear when the tangle $F$ goes over and under
the arc of axis, in steps between (1) and (2), (3) and (4) in the figure.
The 
quandle colors assigned to $F$ changes when $F$ goes under the axis, 
between (3) and (4), and every 
quandle
element $b$ assigned to an arc in $F$ 
changes to 
the quandle element 
$b*a$ if the arc of axis is 
colored by $a$.

\begin{figure}[h]
\begin{center}
\mbox{
\epsfxsize=5in
\epsfbox{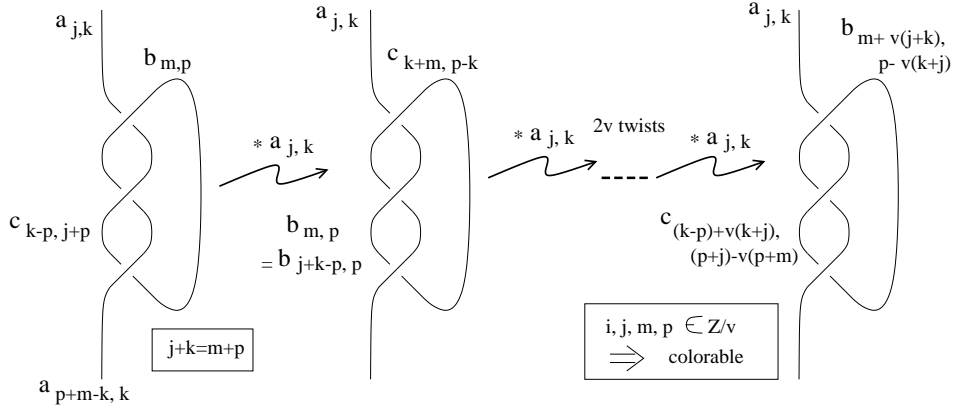} 
}
\end{center}
\caption{The $2v$-twist spun trefoil} 
\label{ntwists} 
\end{figure}

First we consider even twist spun trefoils. 
If we color the arcs of the trefoil as in the top left of
 Fig.~\ref{ntwists} 
 with color 
$a_{j,k}$
on the main arc 
of the axis 
of rotation, 
and color 
$b_{m,p}$ 
on the right, then such a coloring extends 
to the entire trefoil if $j+k=m+p$. 
After each pair of twists the indices $a$, $b$, and $c$ return to
 the arcs of the diagram, but the subscripts $m,p,k-p$ and $j+p$ are
 incremented to $k+j+m$, $p-k-j$, $k+m$, and $p-k$ respectively. 
After $2v$ full twists, the colors 
on the $b$-arc and the $c$-arc become 
$b_{m+v(j+k),\ p-v(k+j)}$ and $c_{(k-p)+v(k+j), \ (p+j)-v(p+m)}$
as indicated in the figure. 
(The subscripts in the figure 
are subjected to the identity $j+k=m+p$ 
 to obtain this result.)   
Thus the extension colors the $2v$-twist spun
 trefoil  if  $v \equiv 0$.

Next we consider a coloring of the trefoil by $QS_4$. 
We label the elements of $(QS)_4(u)$ as follows:
$$[0,j,k,\ell]
=\left( \begin{array}{cccc}
1&0&0&0\\ 
0&0&x^j&0\\
0&0&0&x^k \\ 
0&x^\ell&0&0
\end{array}\right), \quad
[i,0,k,\ell]
=\left( \begin{array}{cccc}
0&0&0&x^i\\ 
0&1&0&0\\
x^k&0&0&0 \\ 
0&0&x^\ell&0
\end{array}\right) $$

$$[i,j,0,\ell]
=\left( \begin{array}{cccc}
0&x^i&0&0\\ 
0&0&0&x^j\\
0&0&1&0 \\ 
x^\ell&0&0&0
\end{array}\right), \quad
[i,j,k,0]
=\left( \begin{array}{cccc}
0&0&x^i&0\\ 
x^j&0&0&0\\
0&x^k&0&0 \\ 
0&0&0&1
\end{array}\right) $$

\begin{figure}[h]
\begin{center}
\mbox{
\epsfxsize=3in
\epsfbox{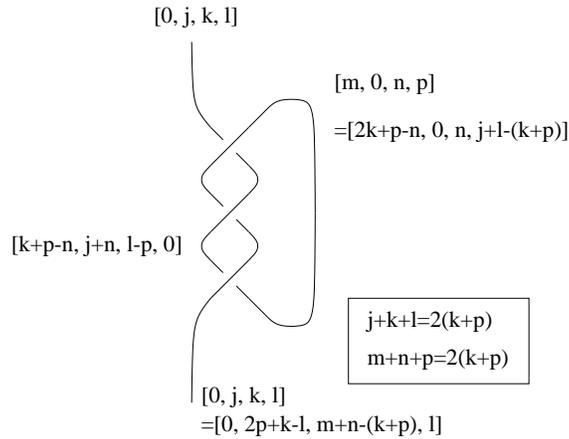} 
}
\end{center}
\caption{Coloring the trefoil with $\widetilde{(QS)}_4$}
\label{trefs4} 
\end{figure}

In the illustration of Fig.~\ref{trefs4}, an extension to a coloring
 by $(QS)_4(u)$ is indicated. The main arc is colored 
$[0,j,k,\ell]$ and the right hand arc is colored $[m,0,n,p].$
These colors induce the color $[k-n+p,j+n,\ell-p,0]$ on the remaining arc of the diagram. A sufficient condition for this to be a coloring 
of trefoil
is that $j+k+\ell = m+n+p=2(k+p).$

We follow the coloring of the right hand arc (colored $[m,0,n,p]$)
for  6 twists.
\begin{eqnarray*} 
[m,\ 0,\ n,\ p] &\stackrel{*[0,j,k,\ell]}{\rightarrow}&
[\ell+m,\ k-\ell+p,\ 0,\ n-k] \\
&\stackrel{*[0,j,k,\ell]}{\rightarrow}&
[j+\ell+m,\ n-k-\ell,\ k+p-j,\ 0] \\
&\stackrel{*[0,j,k,\ell]}{\rightarrow}&
[j+k+\ell+m,\ 0,\ -j-k-\ell+n,\ p] \\
& \stackrel{*[0,j,k,\ell]}{\rightarrow}&
[j+k+2\ell+m,\ k-\ell+p,\ 0,\ n-j-2k-\ell] \\
&\stackrel{*[0,j,k,\ell]}{\rightarrow}& 
[2j+k+2\ell+m,\ -j-2k-2\ell+n,\ -j+k+p,\ 0]\\
&\stackrel{*[0,j,k,\ell]}{\rightarrow}& 
[2(j+k+l)+m,\ 0,\ n-2(j+k+\ell),\ p]\end{eqnarray*}

 Observe that the 3-twist spun trefoil 
colors 
non-trivially 
with $QS_4(2)$, since 
$2(k+p) =0 =(j+k +\ell)$ in this case. 
The result follows by induction. 

A similar calculation applies 
 to the figure $8$ knot and $R_5$.
 We leave the details to the reader. \qed

\section{Variations of cocycle knot invariants} \label{invsec}

The following variations of cocycle knot invariants for 
classical knots have been considered.

\begin{itemize}

\item
For a link $L=K_1 \cup \cdots \cup K_n$, let ${\cal T}_i$, 
$i=1, \ldots, n$, be the set of crossings at which the under-arcs
belong to the component $K_i$. Then it was observed \cite{CENS} 
that $\vec{\Phi}(K) = ( \sum_{\cal C} \prod_{ \tau \in {\cal T}_i }  B(\tau, {\cal C} ) )_{i=1}^n $ is a link invariant, strictly stronger than 
the single state-sum.

\item
Lopes~\cite{Lopes}  
observed that the family
$\{ \prod_{\tau} B(\tau, {\cal C} ) \}_{\cal C \in \mbox{Col}} $ 
is a knot invariant, without taking summation.
Here, $\mbox{Col} $ denotes the set of colorings. 
In particular, infinite quandles can be used for coloring
in this case.
He also defined for links  $L=K_1 \cup \cdots \cup K_n$
the vector version 
$( \{ \prod_{\tau \in {\cal T}_i} 
B(\tau, {\cal C} ) \}_{\cal C \in \mbox{Col}} )_{i=1}^n $. 
\end{itemize}

Now we combine these variations to define the following 
generalized cocycle invariant.

\begin{definition}{\rm   
Let $X$ be a  quandle, $\phi \in {\Z}_{\rm Q}^2(X;A)$, 
where $A$ is an abelian group, ${\cal C}$ a coloring of $L$ 
by $X$, and $B(\tau ,{\cal C})$ the Boltzmann weight at a crossing 
$\tau$ for a coloring ${\cal C}$.
Let $L=K_1\cup \cdots \cup K_r$ be a link and ${\cal T}_i$, 
$i=1,\ldots, r$, be the set of crossings of $L$ such that 
the under-arcs belong to $K_i$. Define
$$\vec{\Psi}(L)=\left\{ \left( \prod_{\tau \in {\cal T}_1}B(\tau ,{\cal C}), 
\ldots, \prod_{\tau \in {\cal T}_r}B(\tau ,{\cal C})
\right) \right\}_{\cal C \in \mbox{Col} }. $$ 
}\end{definition}

This version of a family of vectors 
is potentially stronger than Lopes's version of a vector of families.
For example, the two distinct families of vectors $\{ (1,t), (t,1) \}$ 
and $\{ (1,1), (t,t) \}$ give rise to the same vector of 
families $( \{ 1, t \}, \{ 1, t \} )$. 
As examples, we evaluate the invariants for Whitehead link and Borromean rings,
using extension cocycles constructed  in 
Section~\ref{extsec}.
We use the coefficient group $A=\Z_q=\{ t^n| n=0, 1, \ldots, q-1\}$
 for a positive integer $q$.

\begin{figure}[h]
\begin{center}
\mbox{
\epsfxsize=1.7in
\epsfbox{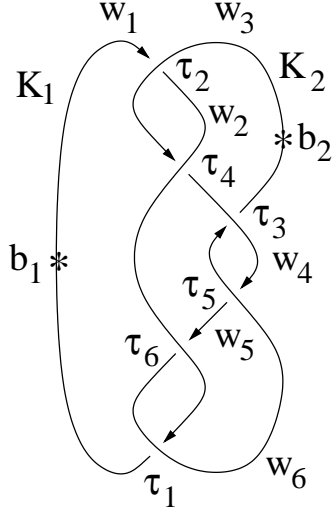}
}
\end{center}
\caption{The Whitehead link} 
\label{fig:whitehead}
\end{figure}

\begin{example} \label{wlinkprop} {\rm
Let $X=W_m={\Z}_{q} [T, T^{-1}] / (1-T)^m$ or
$X=U_m={\Z}_{q^m} [T, T^{-1}] / (T-1+q)$,
and $L$ the Whitehead link. Then the generalized cocycle invariant is
$$\vec{\Psi}(L)= \left\{ \begin{array}{ll}
                      \{ \underbrace{(1,1),\ldots,(1,1)}_{q^{2m}\ copies} 
\}    
& \mbox{for $m=1, 2$,} \\
\{ \underbrace{(t^{n},t^{-n}),\ldots,(t^{n},t^{-n})}_{q^{m+2}\        
copies}\}_{n \in \{0,1,\ldots,q-1\}} & \mbox{for $m\ge3$.}                     
                    \end{array}
            \right.$$
Consequently, 
$$\vec{\Phi}(L)= \left\{ \begin{array}{ll}
                      (q^{2m},q^{2m}) & \mbox{for $m=1, 2$,} \\
(q^{m+2}(t^{q-1}+\cdots+t+1), q^{m+2}(t^{q-1}+\cdots+t+1)) & \mbox{for 
$m\ge3$.} 
                    \end{array}
            \right.$$  
} \end{example}      
{\it Proof.\/}  
Let $X={\Z}_{q} [T, T^{-1}] / (1-T)^m$.
The case for   $X={\Z}_{q^m} [T, T^{-1}] / (T-1+q)$ is similar.
Pick base points $b_1$ and $b_2$ on the components $K_1$ and $K_2$, 
respectively, of the Whitehead link $L=K_1 \cup K_2$ as depicted 
in Fig.~\ref{fig:whitehead}, and trace each component in the given 
orientation of the link. The colors (elements of $X$)
assigned to the arcs that appear in 
this order are $w_1, w_2$ for $K_1$, and $w_3, \ldots, w_6$ for
$K_2$ as depicted. The crossing at the initial point of the arc
colored by $w_i$ is defined to be $\tau_i$.
First we determine the set of colorings: 
For $m\geq 3$, 
for two elements $w_1$, $w_3 \in X$ 
assigned to the top two arcs
of the Whitehead link $L$ as shown in Fig.~\ref{fig:whitehead}, 
there is a coloring of $L$ by $X$ which restricts to the given 
 $w_1$, $w_3 $ 
if and only if 
$$w_3-w_1 \equiv 0 \pmod {(1-T)^{m-3}}\ 
 \mbox{\rm for} \ \    w_1, w_3 \in {\Z}_{q} [T, T^{-1}] / (1-T)^m . $$
 For $m=1,2$, there is such a coloring for any 
$w_1,w_3 \in X$. 
This can be computed as follows.

Represent the elements of $X=\Z_{q} [T, T^{-1}] / (1-T)^m$
by $ a= a_{m-1}(1-T)^{m-1} + \cdots + a_1(1-T) + a_0 $, where
$ a_j \in {\Z}_q $. 
Note that $(1-(1-T))(1+(1-T)+\cdots+(1-T)^{m-1})=1$ in 
$X$, so 
$T^{-1} = 1 + (1-T) + \cdots + (1-T)^{m-1}$. Note also that 
$a\ \bar{*}\ b=T^{-1}a+(1-T^{-1})b$. 

We have the following calculations for each arc:
\begin{eqnarray*}
w_2& =& w_1*w_3 = w_1 + (1-T)(w_3-w_1)\\
w_4 &= & w_3*w_2 = w_3 + (1-T)(w_2-w_3) 
= (w_3-w_1)(1-T)^2 -(w_3-w_1) (1-T) + w_3 \\
w_6 &=& w_3\ \bar{*} \ w_4
 = T^{-1}w_3 + (1-T^{-1})w_4 = (w_3-w_1)(1-T)^2+w_3 \\
w_5 &=& w_4*w_6 = w_4 + (1-T)(w_6-w_4) = 2(w_3-w_1)(1-T)^2-(w_3-w_1)(1-T)+w_3.
\end{eqnarray*}
These relations are obtained using the top four crossings 
($\tau_2, \tau_4, \tau_3$, and $\tau_5$, respectively). 
The  bottom two crossings 
($\tau_6$ and $\tau_1$) 
of the link give rise to 
the next two 
relations.
The first relation is $w_6*w_2 = w_5$ 
 for the second bottom crossing,
giving 
$(w_1-w_3)(1-T)^3 \equiv 0 \pmod{(1-T)^m}$. 
The second relation that corresponds to 
the bottom  crossing is $w_1*w_6 = w_2$ 
 giving 
 $(w_3-w_1) (1-T)^3\equiv 0 \pmod{(1-T)^m}$, as claimed above.

Now we determine the contribution to the invariant for each coloring.
Recall that  
$\phi (w_1,w_3)=[s(w_1)*s(w_3)-s(w_1*w_3)]/(1-T)^m$. 
Since 
$(w_3-w_1)(1-T)^3 \equiv 0 \pmod{(1-T)^m}$, 
we see that the contribution is
\begin{eqnarray*}
\lefteqn{\phi(w_1,w_3)-\phi(w_1,w_6)} \\
&=&
[s(w_1)*s(w_3)-s(w_1*w_3)]/(1-T)^m \\ & & -
[s(w_1)*s(w_3)+(w_3-w_1)(1-T)^3 -s(w_1*w_3)]/(1-T)^m \\
&=&
-(w_3-w_1)(1-T)^3/(1-T)^m \pmod{q}, 
\end{eqnarray*}
for the first component, and for the second component, computations 
show that  
$$\phi(w_3,w_2)-\phi(w_6,w_2)+\phi(w_6,w_4)+\phi(w_4,w_6) =
(w_3-w_1)(1-T)^3/(1-T)^m \pmod{q}. 
$$
 
For $m=1,2$ the contributions for the first and the second component
are both $0$, and 
we have $q^m$ choices for both $w_1$ and $w_3$, 
therefore $\vec{\Psi}(L)=( \underbrace{(1,1),\ldots, 
(1,1)}_{q^{2m}\ copies} )$ and 
$\vec{\Phi}(L)=(q^{2m},q^{2m})$.  

For $m \ge 3$, if $w_1$ and $w_3$ 
color $L$, then 
$(w_3-w_1)(1-T)^3$ 
is $0$
as an element of $X$,
so that $w_3-w_1$  
is uniquely written as 
$w_3-w_1=k(1-T)^{m-3}$, 
where $k=k_0+k_1(1-T)+k_2(1-T)^2$, and 
$k_0, k_1, k_2 \in \{0, 1,\ldots, q-1\}$.  
Then 
$$
(w_3-w_1)(1-T)^3 = k(1-T)^m 
  = (k_0+k_1(1-T)+k_2(1-T)^2)(1-T)^m 
  = k_0(1-T)^m\ \in E.$$ 
Thus the contribution to the 
invariant for the first and second components are 
$t^{-k_0}$ and  $t^{k_0}$, respectively. 

To find the number of colorings conrtibuting to $t^{-k_0}$ and  
$t^{k_0}$, fix $k_0$.  
We have $q^m$ choices for $w_1$ 
 and $q^2$ choices for $k$. 
Then $w_3$ 
is uniquely determined by $w_3=w_1+k(1-T)^{m-3}$. 
 In total, the contribution is $q^m q^2=q^{m+2}$ 
for each $t^{-k_0}$ and $t^{k_0}$. Setting $n=-k_0$ we  obtain the result.
\qed

\begin{figure}[h]
\begin{center}
\mbox{
\epsfxsize=2in
\epsfbox{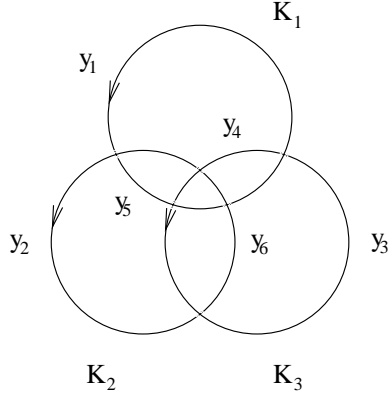}
}
\end{center}
\caption{Borromean Rings} 
\label{fig:rings}
\end{figure}

\begin{example} \label{ringsprop} {\rm 
Let 
$X={\Z}_{q} [T, T^{-1}] / (1-T)^m$ or
 $X={\Z}_{q^m} [T, T^{-1}] / (T-1+q)$,
and $L$ the Borromean rings. Then the generalized cocycle invariant is
$$\vec{\Psi}(L)= \left\{ \begin{array}{ll} 
                 \{ \underbrace{(1,1,1),\ldots,(1,1,1)}_{q^{3m}\ copies} \}
& \mbox{for $m=1$,} \\
\{ \underbrace{(t^{-k_0},t^{-\ell_0},t^{k_0+\ell_0}),\ldots,                       
(t^{-k_0},t^{-\ell_0},t^{k_0+\ell_0})}_{q^{m+2}\ copies}\}_{k_0,\ell_0 
\in                
\{0,1,\ldots,q-1\}} & \mbox{for $m\ge2$.}
                    \end{array}                                                
            \right.$$
Consequently,
$$\vec{\Phi}(L)= \left\{ \begin{array}{ll}
                      (q^{3m},q^{3m},q^{3m}) & \mbox{for $m=1$,} \\
(q^{m+2}(t^{q-1}+\cdots+t+1),q^{m+2}(t^{q-1}+\cdots+t+1),
q^{m+2}(t^{q-1}+\cdots+t+1)) & \mbox{for $m\ge2$.} 
                    \end{array}
            \right.$$
} \end{example}
{\it Proof.\/}  
Let $X={\Z}_{q} [T, T^{-1}] / (1-T)^m$ and let $L$ be the Borromean
rings as depicted in Fig.~\ref{fig:rings}. The case for  
 $X={\Z}_{q^m} [T, T^{-1}] / (T-1+q)$ is similar.
Calculations are similar to the preceding example 
and we give a sketch.
First we determine the set of colorings: For three elements 
$y_1,y_2,y_3 \in X$  assigned to each outer arc in the diagram of $L$, 
there is a coloring of $L$ by $X$ which restricts to the given 
$y_1,y_2,y_3$ if and only if
$$(y_2-y_3) (1-T)^2 \equiv 0, \  \
(y_1-y_2)(1-T)^2 \equiv 0 \ \ \mbox{and} \ \ 
(y_3-y_1) (1-T)^2 \equiv 0 . $$
The outer three crossings are used to describe $y_4, y_5, y_6$ 
in terms of the rest, and the inner three crossings give the above relations.

Contributions to the invariant are computed as follows.  
The contribution for the first component of $L$ colored by $y_1$ is
$\phi(y_1,y_2)-\phi(y_1,y_2*y_3)=-(y_3-y_2)(1-T)^2/(1-T)^m \pmod{q}$. 
For $m=1$, the contribution is trivial, and  the total number of colorings
is $q^{3m}$.
For $m\ge2$, $(y_3-y_2)(1-T)^2$ 
is divisible by $(1-T)^m$, so $y_3-y_2$ is uniquely 
written      
as $y_3-y_2=k(1-T)^{m-2}$, where $k=k_0+k_1(1-T)$ and 
$k_0, k_1 \in \{0,1,\ldots,q-1\}$.
So $$
(y_3-y_2) (1-T)^2= k(1-T)^m
  = (k_0+k_1q)(1-T)^m
  = k_0(1-T)^m ,
$$
and the first component contributes $t^{-k_0}$ to the invariant.
For the second component of $L$ colored by $y_2$, similar calculations
as above give 
the contribution $\phi(y_2,y_3)-\phi(y_2,y_3*y_1)=-(y_1-y_3)(1-T)^2$,
which  is divisible by $(1-T)^m$ so $y_1-y_3=(\ell_0+\ell_1(1-T))(1-T)^{m-2}$
and therefore $-(y_1-y_3)(1-T)^2=-\ell_0 (1-T)^m$.
The we obtain
$y_2-y_1=-[(k_0+\ell_0)+(k_1+\ell_1)(1-T)](1-T)^{m-2}$, 
so that the third component contributes  
 $t^{k_0+\ell_0}$.
Finally, the contribution to the invariant is the vector
$(t^{-k_0}, t^{-\ell_0}, t^{k_0+\ell_0})$, where the entries
correspond to the components $K_1, K_2, K_3$, respectively.
The result follows. \qed

In the above examples, we see that the cocycle invariant is non-trivial
when the given link is colored by 
$X={\Z}_{q} [T, T^{-1}] / (1-T)^m$ but not by 
$E={\Z}_{q} [T, T^{-1}] / (1-T)^{m+1}$,
and the descrepancy in extending the coloring contributes to
the invariant. This is the case in general, as
proved in \cite{CENS} for the knot case. 
We rephrase the theorem in our situation and include a 
similar proof for reader's convenience.

Let $K$ be a classical or virtual  knot or link. 
Let ${\cal C}$ be a coloring of $K$ by $X$. Let $E$ be
a quandle with a surjective homomorphism 
$p: E \rightarrow X$. 
If there is a coloring ${\cal C}'$ of $K$ by $E$ such that for 
every arc $a$ of $K$, it holds that $p ({\cal C}'(a))={\cal C}(a)$, 
then  ${\cal C}'$ is called an  {\it extension} of  ${\cal C}$.

\begin{theorem} \label{linkcolorthm}
Let 
 $$\vec{\Psi}(L)=\left\{ \left( \prod_{\tau \in {\cal T}_1}B(\tau ,{\cal C}), 
\ldots, \prod_{\tau \in {\cal T}_r}B(\tau ,{\cal C})
\right) \right\}_{\cal C \in \mbox{\rm Col} }$$  
be the generalized 
cocycle invariant of a link $L=K_1 \cup \cdots \cup K_n$
with a quandle $X$ and a cocycle $\phi \in Z^2_{\rm Q}(X;A)$ for 
an abelian group $A$. 
Then  
$\left( \prod_{\tau \in {\cal T}_1}B(\tau ,{\cal C}), 
\ldots, \prod_{\tau \in {\cal T}_r}B(\tau ,{\cal C})
\right)$
 is 
 a vector with every entry $1$ 
for a 
coloring ${\cal C}$ if and 
only if 
the coloring ${\cal C}$ 
extends to a coloring of $L$ by $E(X, A, \phi)$. 
\end{theorem}
{\it Proof.\/}
 Let ${\cal C}$ be a coloring whose contribution to 
$\vec{\Psi}(L) $ is  $(1, \ldots, 1)$. 
Fix this coloring in what follows.
Pick a base point $b_0$ on a component $K_i$ of $L$.
Let $x \in X$ be the color on the arc $\alpha_0$ containing $b_0$.
Let $\alpha_i$, $i=1, \ldots, n$, be the set of arcs that appear 
in this order when the diagram $K$ is traced
 in the given orientation of $K_i$, 
starting from $b_0$.
Pick an element $a \in A$ and give a color 
$(a, x)$ on  $\alpha_0$,
so that we define a coloring ${\cal C}'$ by $E$ on $\alpha_0$
by ${\cal C}'(\alpha_0)=(a, x)\in E$. 
We try to extend it to the entire
diagram by traveling the diagram from $b_0$ along the arcs $\alpha_i$,
 $i=1, \ldots, n$, in this order, by induction. 

Suppose  ${\cal C}'(\alpha_i)$ is defined for $0 \leq i < k$. 
Define  ${\cal C}'(\alpha_{k+1})$ as follows.
Suppose that the crossing $\tau_k$ separating $\alpha_k$ and $\alpha_{k+1}$
is positive, and the over-arc at  $\tau_k$ is $\alpha_j$. 
Let  ${\cal C}'(\alpha_{k})=(a, x)$ and ${\cal C}(\alpha_j)=y \in X$.
Then we have ${\cal C}(\alpha_{k+1})=x*y \in X$. 
Define  ${\cal C}'(\alpha_{k+1})=(a  \phi(x, y), x*y)$ in this case.

 Suppose that the crossing $\tau_k$ is negative.
Let  ${\cal C}'(\alpha_{k})=(a, x)$ and ${\cal C}(\alpha_j)=y \in X$.
Then if  ${\cal C}(\alpha_{k+1})=z$, then we have $z*y=x$.
Define  ${\cal C}'(\alpha_{k+1})=(a  \phi(z, y)^{-1}, z)$ in this case.

Define ${\cal C}'(\alpha_{i})$ inductively for all $i=0, \ldots, n$.
Regard $\alpha_0$ as $\alpha_{n+1}$, and repeat the 
above 
construction at the last crossing $\tau_n$ to come back to $\alpha_0$.
By the construction we have 
 ${\cal C}'(\alpha_{n+1})=(a  \prod_{\tau} B(\tau, {\cal C}),
{\cal C}(\alpha_0) )$, where  $\prod_{\tau} B(\tau, {\cal C})$
is the state-sum contribution (the product of Boltzmann weights over all crossings) of  ${\cal C}$.
This   contribution is  equal to $1$ by the assumption that 
$\prod_{\tau} B(\tau, {\cal C})=1$,
and we have a well-defined coloring ${\cal C}'$. 
Hence this color extends to $E(X,A,\phi)$.

Conversely, if a coloring
${\cal C}$ by $X$ extends to a coloring by $E(X,A,\phi)$, 
then from the above argument, we have that  
$(a, x)=(a \prod_{\tau} B(\tau, {\cal C}), x )$,
if $(a, x)$ is the color on the base point $b_0$.
Hence $ \prod_{\tau} B(\tau, {\cal C})=1$.
$\qed$

\section{Cocycle invariants and Alexander 
 matrices} \label{LXsec}

In  this section we point out 
relations of the 
cocycle invariants to Alexander 
matrices.
We examine closely Example~\ref{wlinkprop} 
given in Section~\ref{invsec} 
from this new point of view.

Let {\em $B_{D_L}=\sum_{i=1}^{n}B_i$} be an $(n \times n)$-matrix 
where $B_i$ is the $(n \times n)$-matrix corresponding to each 
crossing point $\tau _i$ 
such that $(k_i,i)$ entry is 
$T^{\epsilon _i}$, $(\ell_i,i)$ entry is $1-T^{\epsilon _i}$ 
and otherwise is 0. Here $\epsilon _i$ means the sign of the 
crossing point $\tau_i$.
Set $A_{D_L}=B_{D_L}-E_n$, 
where $E_n$ denotes the $n$-dimensional identity matrix, then 
from the definitions it follows that    
 $A_{D_L}$ is an Alexander matrix.
Recall that a coloring is a function ${\cal C} : R \rightarrow X$, 
where $R$ is the set of over-arcs in the diagram and $X$ is a fixed 
Alexander quandle $\Lambda/J$ for an ideal $J$. 
A coloring which assigns $y_i$ to an arc $a_i$ (${\cal C}(a_i)=y_i$) is 
represented by the vector $\vec{y}=(y_1, \ldots, y_n)$ satisfying 
$\vec{y}A_{D_L}^{(X)}=\vec{0}$.   These 
descriptions are given in \cite{Inoue} to prove Theorem~\ref{Inouethm}.

\begin{proposition} \label{numofquandlehomos}       
Let $L=K_1 \cup \cdots \cup K_r$ be a link and $X = \Lambda_q /J$     
be an Alexander quandle. Suppose $E= \Lambda_{q'} / J'$ is an abelian 
extension of $X$, where $q, q'$ are positive integers. Let 
$A^{(X)}_{D_L}$ (respectively $A^{(E)}_{D_L}$) be the matrix 
$A_{D_L}$ regarded as a matrix over $X$ (respectively over $E$).
Then a coloring 
$\vec{y}$ of $L$ by $X$ 
contributes 
a non-trivial value to the invariant $\vec{\Psi}(L)$
if and only if $\vec{y}A^{(X)}_{D_L} =\vec{0}$ and 
$s(\vec{y})A^{(E)}_{D_L}=\vec{x}\ne \vec{0}$,
where $s: X \rightarrow E$ is the natural section.   
\end{proposition}
{\it Proof.\/} 
Let $\psi : (\Lambda _q/J)^n \rightarrow (\Lambda _q/J)^n$ be the map 
which takes a row vector $\vec{y}$ to $\vec{y}A_{D_L}$.
By Inoue's description given above, 
 the set of all quandle colorings is equal to 
ker$A^{(X)}_{D_L}$.
If $ \vec{y}A^{(X)}_{D_L}=\vec{0}$ and
 $s(\vec{y})A^{(E)}_{D_L}=\vec{x}\neq 0$, 
then by Theorem~\ref{linkcolorthm}
we obtain that  $\vec{\Psi}(L)$ is non-trivial.
$\qed$\\

\begin{figure}
\begin{center}
\mbox{
\epsfxsize=1.5in
\epsfbox{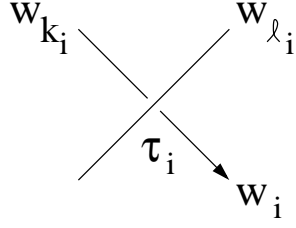} 
}
\end{center}
\caption{ Labeling a crossing }
\label{colorcrossing} 
\end{figure}

Next, we compute the non-trivial contributions using Alexander matrices,
for extensions discussed in Section~\ref{extsec}.
Let $X = W_m=\Lambda_q /(1-T)^m$
or $X=U_m= \Lambda_{q^m}/(T-1+q)$,
and   $E= W_{m+1}$ or $E=U_{m+1}$  be their abelian extensions,
respectively.
For this purpose, we fix the following convention in 
numbering crossings and arcs of a given diagram.

 Let $L=K_1 \cup \cdots \cup K_r$ be a link with $n$ crossings. 
Pick a base point $b_i$ on $K_i$, for $i=1, \ldots, r$. 
Let $a_1, \ldots, a_{i_1}$ be the arcs of $K_1$ such that $a_1$ contains 
$b_1$ and they appear in this order when one traces $K_1$ in the given 
orientation of $K_1$ starting from $b_1$. Then let $a_{i_1+1}$ be the 
arc of $K_2$ containing $b_2$ and $a_{i_1+2}, \ldots, a_{i_2}$ be the 
arcs of $K_2$ similarly defined from the given orientation. Repeat this 
process for the remaining components to obtain the arcs 
$a_1, \ldots, a_{i_1}, a_{i_1+1}, \ldots, a_{i_2}, a_{i_2+1},  
\ldots, a_{i_{r-1}+1}, \ldots, a_{i_r}=a_n$. 
Let ${\cal C} : R \rightarrow X$ be a coloring of $L$ by $X$. 
Let $w_i={\cal C}(a_i)$ and $\tau_i$ be the crossing such that the 
outcoming under-arc is 
$a_i$ 
for $i=1, \ldots, n$ (see Fig.~\ref{colorcrossing}).
This convention is used in Fig.~\ref{fig:whitehead}.

 Let
$s: X \rightarrow E$ 
be  the 
 section defined in Section~\ref{extsec} respectively by 
\begin{eqnarray*}
s\left( \sum_{j=0}^{m-1} A_j (1-T)^j \quad \mbox{mod}\ (1-T)^{m} \right)
 &= & \sum_{j=0}^{m-1} A_j (1-T)^j \quad  \mbox{mod}\ (1-T)^{m+1}
\quad \mbox{\rm for $W_m$, and} \\
s \left( \sum_{j=0}^{m-1} X_j q^j \right) &= & 0 \cdot q^{m}+ 
\sum_{j=0}^{m-1} X_j q^j \quad \mbox{\rm for $U_m$.}
\end{eqnarray*} 
For the following theorem, let
 $\vec{\Psi}(L)$ be the generalized cocycle invariant
defined with the cocycle $\phi \in Z^2_{\rm Q}(X;{\Z}_q)$ 
corresponding to the extension $p: E \rightarrow X$
specified above.

\begin{proposition} \label{contributionthm}
Let $A_{D_L}$ be the Alexander matrix 
obtained from $D_L$ with the above choice of order of $w_i$ 
and $\tau_i$. 

A given coloring represented by a vector
$\vec{w}$ contributes 
a non-trivial vector to the invariant $\vec{\Psi}(L)$
if and only if $\vec{w}A^{(X)}_{D_L} =\vec{0}$ and 
$s(\vec{w})A^{(E)}_{D_L}=\vec{z}\ne \vec{0}$. This contribution is
\begin{eqnarray*}
(t^{\sum_{j=1}^{i_1}\eta(\tau_j) z_j/(1-T)^m}, 
\ldots ,t^{\sum_{j=i_{r-1}+1}^{i_r}\eta(\tau_j) z_j/(1-T)^m}) 
 & & \mbox{\rm for $X=W_m$, and } \\
(t^{\sum_{j=1}^{i_1}\eta(\tau_j) z_j/q^m},  
\ldots ,t^{\sum_{j=i_{r-1}+1}^{i_r}\eta(\tau_j) z_j/q^m}) 
 & & \mbox{\rm for $X=U_m$, respectively}
\end{eqnarray*}
where $\eta(\tau)=1$ for a positive crossing $\tau$ and
$\eta(\tau)=T$ for a negative crossing $\tau$. 
\end{proposition}
{\it Proof.\/} 
We consider the case $X=W_m$, as the other case is similar.
Let $\psi : (\Lambda _q/(1-T)^m)^n \rightarrow (\Lambda _q/(1-T)^m)^n$ 
be the map which takes a row vector $\vec{w}$ to $\vec{w}A_{D_L}^{(X)}$.
 Assume that $\vec{w}A_{D_L}^{(X)}=\vec{0}$ and 
$s(\vec{w})A_{D_L}^{(E)}=\vec{z} \ne \vec{0}$. 
The contribution to the invariant at a 
positive 
crossing $\tau_i$ is given by 
\begin{eqnarray*}
\phi(w_{k_i},w_{\ell_i}) 
&=&[s(w_{k_i})*s(w_{\ell_i})-s(w_{k_i}*w_{\ell_i})]/(1-T)^m \\
   &=& [s(w_{k_i})*s(w_{\ell_i})-s(w_i)]/(1-T)^m ,
\end{eqnarray*}
where $w_{\ell_i}$ is the color on the over-arc at the crossing $\tau_i$, 
and $w_{k_i}$ is the color on the incoming under-arc at $\tau_i$ 
if $\tau_i$ is positive 
(see Fig.~\ref{colorcrossing}).
Since $\vec{w}$ is in the kernel, 
$w_{k_i}*w_{\ell_i} - w_i= Tw_{k_i}+(1-T)w_{\ell_i} - w_i=0\ \mbox{mod}(1-T)^m$
and we have $[s(w_{k_i})*s(w_{\ell_i})-s(w_i)]/(1-T)^m =z_i/(1-T)^m$.

Suppose $\tau_i$ is negative. Then the contribution is
\begin{eqnarray*}
 - \phi(w_{i},w_{\ell_i}) 
&=&- [s(w_{i})*s(w_{\ell_i})-s(w_{i}*w_{\ell_i})]/(1-T)^m \\
   &=& -[s(w_{i})*s(w_{\ell_i})-s(w_{k_i})]/(1-T)^m \\
 &=& -[ T w_{i} + (1-T) w_{\ell_i} - w_{k_i} ] /(1-T)^m.
\end{eqnarray*}
On the other hand, 
$$ z_i= T^{-1}  w_{k_i} + (1-T^{-1} ) w_{\ell_i} -  w_{i}
= -T^{-1} [  T w_{i} + (1-T) w_{\ell_i} - w_{k_i} ]$$
so that the contribution is $T z_i$ in this case. 
Hence the total contribution of the invariant for the component $K_r$ is
$$ t^{\sum_{j=i_{r-1}+1}^{i_r} \eta(\tau_j)z_j/(1-T)^m}, $$ 
where $\{z_1, \ldots, z_{i_r}\} \in K_r$.
$\qed$

\begin{example}{\rm
\begin{sloppypar}
We consider  the Whitehead link $L=K_1 \cup K_2$ 
depicted in Fig.~\ref{fig:whitehead}.
Let $X=W_m$ and $E=W_{m+1}$. 
Use the letters $w_i$, 
($i=1, \ldots, 6$) as depicted in the figure as 
colors assigned to the arcs, as well as generators for the Alexander matrix.
 Then the 
Alexander matrix $A_{D_L}=B_{D_L}-E_n$
with respect to  the columns corresponding  to 
$(\tau_1,\ldots,\tau_6)$
and rows  corresponding  to $(w_1, \ldots, w_6)$ is
given by 
\end{sloppypar}
$$ A_{D_L}=\left(
\begin{array} { c c c c c c }
-1    & T   & 0 & 0    & 0 & 0 \\
T^{-1} & -1  & 0 & 1-T & 0 & 1-T^{-1} \\
0      & 1-T & -1& T   & 0 & 0 \\
0     & 0    &1-T& -1  & T & 0  \\
0     & 0    & 0 & 0   & -1 & T^{-1}  \\
1-T^{-1} & 0 & T & 0   & 1-T & -1
\end{array} 
\right). $$
After some row and column permutations we obtain 
$$ A_0=\left(
\begin{array} { c c c c c c }
-1 & 1-T & 0 & 0 & 1-T^{-1} & T^{-1} \\
0 & -1& 1-T & T & 0 & 0 \\
0 & 0 & T & 1-T & -1 & 1-T^{-1}  \\
0 & 0 & 0 & -1 & T^{-1} & 0  \\
T & 0 & 0 & 0 & 0 & -1  \\
1-T & T & -1 & 0 & 0 & 0 
\end{array} 
\right), $$
\sloppypar
\noindent
with respect to 
 the columns corresponding  to 
$(\tau_2, \tau_4, \tau_3, \tau_5, \tau_6, \tau_1)$
and rows  corresponding  to $(w_2,w_4,w_6,w_5,w_1,w_3)$.
This permutation is performed so that
 we can diagonalize the first four rows and 
columns by column reductions
to obtain
\sloppypar
$$ A_1 =\left(
\begin{array} { c c c c c c }
1 & 0 & 0 & 0 & 0 & 0 \\
0 & 1 & 0 & 0 & 0 & 0 \\
0 & 0 & 1 & 0 & 0 & 0 \\
0 & 0 & 0 & 1 & 0 & 0  \\
-T & -T+T^2 & (1-T)^2 & 1-3T+2T^2 & -T^{-1}(1-T)^3 & T^{-1}(1-T)^3 
\\
 -1+T & -1+T-T^2 & -1-(1-T)^2 & -2+3T-2T^2 & T^{-1}(1-T)^3 & -T^{-1}(1-T)^3 
\end{array} 
\right). $$ 
The solution set 
$$(w_2,w_4,w_6,w_5,w_1,w_3)A_1^{(X)}=(0,0,0,0,0,0), $$
is written by 
\begin{eqnarray*}
w_2 &=&Tw_1+(1-T)w_3 \\
w_4 &=& T(1-T)w_1+(T+(1-T)^2)w_3 \\
w_6 &=& -(1-T)^2w_1+(1+(1-T)^2)w_3 \\
w_5 &=& (T(1-T)-(1-T)^2)w_1+(T+2(1-T)^2)w_3 \\
0 &=& (w_3-w_1)T^{-1}(1-T)^3 
\end{eqnarray*} 
where $A_1^{(X)}$ denotes the matrix $A_1$ regarded as a matrix over $X$.
The  set of colorings is represented by vectors in the kernel
of $A_1^{(X)}$. Specifically, the kernel is the set of vectors
$\vec{w}$
with $w_1$ and $w_3$ 
satisfying    $(1-T)^3(w_3-w_1)=0$ 
in $X$
and $w_2,w_4,w_6,w_5$ 
determined accordingly as above. 
This  matches the computations in Example~\ref{wlinkprop}.
The contribution to the invariant is obtained by computing 
\begin{eqnarray*}
\vec{z} &=& s(\vec{w}) A_{D_L}^{(E)} \\
&=& (-T^{-1}(1-T)^3(w_3-w_1),\ 0,\ 0,\ 0,\ 0,\ T^{-1}(1-T)^3(w_3-w_1)).
\end{eqnarray*} 
By Proposition~\ref{contributionthm} the non-trivial contribution to 
$\vec{\Psi}(L)$ is $(t^{\sum_{j=1}^{2}\eta(\tau_j) z_j/(1-T)^3}, 
t^{\sum_{j=3}^{6}\eta(\tau_j) z_j/(1-T)^3})=(t^{-s},t^{s})$ for 
some $s$ $(0\le s \le q-1)$ depending on the value of $w_3-w_1$, 
and this matches Example~\ref{wlinkprop}.
}\end{example}

Finally we observe a relation to the Conway polynomial.
Let $\Delta_L(T)\in {\Z}[T^{-\frac{1}{2}},T^{\frac{1}{2}}]$
be the {\em Conway-normalized Alexander polynomial}~\cite{Lickorish}.
In our case, let $A^{'}_{D_L}$ be the matrix obtained from
$A_{D_L}$ by deleting the $j$th column and $j$th row for some
$j$, $j=1,\ldots,n$,  
 let $f(T)=\mbox{\rm det}(A^{'}_{D_L}) \in {\Z}[T^{1},T^{-1}]$
and $\mu$ and $\nu$ be the maximal and minimal degree of $f$ respectively.
Then $\Delta_L(T)=T^{-\frac{\mu + \nu}{2}}f(T)$.
The {\em Conway polynomial} $\nabla_L(z) \in {\Z}[z]$ is defined by
$\nabla_L(T^{-\frac{1}{2}}-T^{\frac{1}{2}})=\Delta_L(T)$ where
$z=T^{-\frac{1}{2}}-T^{\frac{1}{2}}$.

\begin{proposition}
Let the minimal degree of $\nabla_L(z)$ be denoted by 
\mbox{\rm min-deg}$\nabla_L(z)$, then it  satisfies
\mbox{\rm min-deg}$\nabla_L(z) \ge m$, where 
$m$ is the smallest integer such that the cocycle invariant
defined from the extension of ${\Z}_q[T,T^{-1}]/(1-T)^m$ to
${\Z}_q[T,T^{-1}]/(1-T)^{m+1}$ is non-trivial.
\end{proposition}
{\it Proof.\/}
Assume that $\vec{y}A^{(X)}_{D_L} =\vec{0}$ and 
$s(\vec{y})A^{(E)}_{D_L}=\vec{x} \ne \vec{0}$. Then $\vec{y}$ contributes 
a non-trivial value to the invariant $\vec{\Psi}(L)$ as in 
Proposition~\ref{numofquandlehomos}.
Since $\vec{x} \ne \vec{0}$ there exists $i$, $1\le i \le n$, 
such that $x_i \ne 0$.
Let $j$ be an integer, $1 \le j \le n$, with $j \ne i$.
Let $\vec{x'}$ be the vector $\vec{x}$ with the $x_j$ entry deleted.
Then there exists $\vec{y'}\ne \vec{0}$, where $\vec{y'}$ is the 
vector $\vec{y}$ with the $j$th entry deleted, such that
$\vec{y'}A^{'(X)}_{D_L} =\vec{0}$.
This implies that $\mbox{\rm det}A^{'(X)}_{D_L} = 0 $. 
Hence $\mbox{\rm det}A^{'}_{D_L} \equiv 0 \pmod{(1-T)^m} $, and we have
\mbox{\rm min-deg}$\nabla_L(z) \ge m$.
$\qed$

\end{document}